\documentclass{article}%
\usepackage{graphicx}
\usepackage{amsmath}
\usepackage{amsfonts}
\usepackage{amssymb}%
\newtheorem{theorem}{Theorem}

\newtheorem{corollary}[theorem]{Corollary}

\newtheorem{lemma}[theorem]{Lemma}

\newtheorem{proposition}[theorem]{Proposition}

\begin{document}

\title{A pinching estimate for solutions of the linearized Ricci flow system on $3$-manifolds}
\author{Greg Anderson\\School of Mathematics\\University of Minnesota
\and Bennett Chow\thanks{Research partially supported by NSF grant DMS-9971891.}\\Department of Mathematics\\University of California at San Diego}
\date{September 17, 2001 (Revised May 3, 2002)}
\maketitle

\section{Introduction}

An important component of Hamilton's program for the Ricci flow on compact
$3$-manifolds is the classification of singularities which form under the flow
for certain initial metrics. In particular, Type I singularities, where the
evolving metrics have curvatures whose maximums are inversely proportional to
the time to blow-up, are modelled on the $3$-sphere and the cylinder
$S^{2}\times\mathbb{R}$ and their quotients. On the other hand, Type II
singularities (the complementary case) are much more difficult to understand.
Despite this, it is known from the work of Hamilton that their singularity
models are stationary solutions to the Ricci flow. This uses several
techniques, including Harnack inequalities of Li-Yau-Hamilton type, the strong
maximum principle for systems, dimension reduction, and the study of the
geometry at infinity of noncompact stationary solutions (see\S \S 14-26 of
\cite{H2}.) In terms of Hamilton's program, at least two obstacles remain:
obtaining an injectivity radius estimate for Type II solutions and ruling out
the so-called cigar soliton (the unique complete stationary solution on a
surface with positive curvature) as the dimension reduction of a Type II
singularity model.\footnote{In fact, Hamilton has announced informally that
these are the only two obstacles and that they both would follow from
obtaining a suitable differential matrix Harnack inequality of Li-Yau-Hamilton
type for arbitrary solutions of the Ricci flow on compact $3$-manifolds.}

On the other hand, it is also conjectured by Hamilton that Type II
singularities are not generic. If this conjecture can be proven with some
definition of generic which implies that for any compact $3$-manifold the
Ricci flow with suitable surgeries (see \cite{H5} for how to perform
surgeries) does not form Type II singularities, then there would be no need
for obtaining an injectivity radius estimate for Type II solutions or ruling
out the so-called cigar soliton. Partly for these reasons we are motivated to
study the linearized Ricci flow. Given an initial metric on a $3$-manifold and
the corresponding solution to the Ricci flow, one would like to understand the
behavior of solutions with nearby initial metrics. The linearized Ricci flow
system is the pair of equations (\ref{ricciflow})-(\ref{LichLapEq}) which we
consider below. In this paper we obtain an apriori estimate for arbitrary
solutions to the linearized Ricci flow on compact $3$-manifolds which we hope
may be useful in its study. The inspiration for this estimate comes from the
works of Hamilton (\S 10 of \cite{H1} and \S 24 of \cite{H2}) and Gursky
\cite{G}.

Recall that if $\left(  M^{3},g\left(  t\right)  \right)  $ is a solution to
the Ricci flow on a compact $3$-manifold with positive scalar curvature, then
Hamilton obtained the following parabolic Bochner-type estimate
\begin{equation}
\frac{\partial}{\partial t}\left(  \frac{\left\vert Rc\right\vert ^{2}}{R^{2}%
}\right)  \leq\Delta\left(  \frac{\left\vert Rc\right\vert ^{2}}{R^{2}%
}\right)  +\frac{2}{R}\nabla R\cdot\nabla\left(  \frac{\left\vert
Rc\right\vert ^{2}}{R^{2}}\right)  .\label{bochner1}%
\end{equation}
See Lemmas 10.5 (with $\gamma=2$) and 10.6 in \cite{H1} for the positive Ricci
curvature case, and the equation for $Y$ in the proof of Theorem 24.7 in
\cite{H2} for the more general positive scalar curvature case. A sharpened
form of this estimate is the main estimate in showing that the normalized
Ricci flow evolves a closed $3$-manifold with positive Ricci curvature into a
spherical space form; see \cite{H1}, Theorem 10.1, or \cite{H3}, Theorem 5.3
for a simpler proof. A further extension of this estimate is used to show that
for a Type I\ singularity of the Ricci flow on a closed $3$-manifold not
diffeomorphic to a spherical space form, there exists a sequence of dilations
about points and times approaching the singularity time that limits to a
quotient of the cylinder $S^{2}\times\mathbb{R}$ ; see Theorem 24.7, Corollary
24.8 and Theorem 26.5 in \cite{H2}.

On the other hand, also recall that Gursky (see \cite{G}) proved that if
$\left(  M^{4},g\right)  $ is a closed, oriented, $4$-manifold with positive
scalar curvature such that%
\begin{equation}
\Delta R=-4\left|  W^{+}\right|  ^{2}-2\left|  Rc-\frac{1}{4}Rg\right|
^{2}+\frac{1}{6}R^{2} \label{Gursky's Eq}%
\end{equation}
and $\alpha$ is a self-dual harmonic $2$-form, then%
\begin{equation}
\Delta\left(  \frac{\left|  \alpha\right|  ^{2}}{R^{2}}\right)  +\frac{2}%
{R}\nabla R\cdot\nabla\left(  \frac{\left|  \alpha\right|  ^{2}}{R^{2}%
}\right)  \geq0. \label{bochner2}%
\end{equation}

There is a formal similarity between these two Bochner formulas. Note that
K\"{a}hler-Einstein surfaces satisfy (\ref{Gursky's Eq}). On the other hand,
the Ricci flow is the parabolic version of the equation for Einstein metrics
(the fixed points of the volume normalized Ricci flow are the Einstein
metrics). Furthermore, if $\left(  M^{2k},g\right)  $ is a K\"{a}hler manifold
and $\alpha$ is a $J$-invariant $2$-form, then $h\left(  X,Y\right)
\doteqdot\alpha\left(  X,JY\right)  $ is a $J$-invariant symmetric $2$-tensor
and $\left(  \Delta_{d}\alpha\right)  \left(  X,JY\right)  =\left(  \Delta
_{L}h\right)  \left(  X,Y\right)  ,$ where $\Delta_{L}$ is the Lichnerowicz
laplacian (defined below). Thus if $\alpha$ is also harmonic, then $\Delta
_{L}h=0.$ The parabolic version of this last equation is the Lichnerowicz
laplacian heat equation.

The above considerations partly motivate us to study the analogue of estimates
(\ref{bochner1}) and (\ref{bochner2}) in the context of the Ricci flow
\begin{equation}
\frac{\partial}{\partial t}g_{ij}=-2R_{ij}\label{ricciflow}%
\end{equation}
coupled to the Lichnerowicz laplacian heat equation
\begin{equation}
\frac{\partial}{\partial t}h_{ij}=\left(  \Delta_{L}h\right)  _{ij}%
\doteqdot\Delta h_{ij}+2R_{kijl}h_{kl}-R_{ik}h_{kj}-R_{jk}h_{ki}%
\label{LichLapEq}%
\end{equation}
for a symmetric $2$-tensor $h.$ This is the linearized Ricci flow system and
arises from linearizing the Ricci flow using a version of DeTurck's trick (see
\S \ref{defmot}). A differential Harnack inequality of Li-Yau-Hamilton type,
patterned after Hamilton's trace inequality for the Ricci flow \cite{H4} and
Li-Yau's seminal inequality for the heat equation \cite{LY}, for this coupled
system was proved by Hamilton and one of the authors \cite{CH} and interpreted
geometrically in terms of linearizing the Ricci flow by S.-C. Chu and one of
the authors \cite{CC}. A complex version of this inequality was proven by Ni
and Tam \cite{NT} and applied to the study of the K\"{a}hler-Ricci flow. The
linearized Ricci flow has been studied by Guenther, Isenberg and Knopf
\cite{GIK} at flat solutions from the point of view of maximal regularity
theory. Additional work on differential Harnack inequalities of
Li-Yau-Hamilton type which appear related to the linearized Ricci flow are in
\cite{C} and \cite{CK}.\footnote{In the latter reference, a Li-Yau-Hamilton
inequality for the Ricci flow is proved which generalizes Hamilton's matrix
inequality and has some formal similarities with linear inequalities.}

Analogous to (\ref{bochner1}) and (\ref{bochner2}), we may consider the
quantity $\frac{\left|  h\right|  ^{2}}{R^{2}}$ for solutions to the Ricci
flow on $3$-manifolds with positive scalar curvature. This is a pointwise
measure of the size of $h$ relative to the scalar curvature. More generally,
since $R_{\min}\left(  t\right)  $ is a nondecreasing function for solutions
to the Ricci flow, we may replace $R$ by $R+\rho,$ where $\rho\in
\lbrack0,\infty)$ is chosen so that $R+\rho>0$ at $t=0.\medskip$

\noindent\textbf{Main Theorem}. \emph{Let }$\left(  M^{3},g\left(  t\right)
\right)  $\emph{ be a solution to the Ricci flow on a closed }$3$%
\emph{-manifold on a time interval }$[0,T)$\emph{ with }$T<\infty$\emph{ and
let }$\rho\in\lbrack0,\infty)$\emph{ be such that }$R_{\min}\left(  0\right)
>-\rho.$\emph{ If the pair }$\left(  g,h\right)  $\emph{ is any solution to
the linearized Ricci flow system (\ref{ricciflow})-(\ref{LichLapEq}), then
there exists a constant }$C<\infty$\emph{ such that}%
\begin{align}
\frac{\partial}{\partial t}\left(  \frac{\left|  h\right|  ^{2}}{\left(
R+\rho\right)  ^{2}}\right)   &  \leq\Delta\left(  \frac{\left|  h\right|
^{2}}{\left(  R+\rho\right)  ^{2}}\right)  +\frac{2}{R+\rho}\nabla
R\cdot\nabla\left(  \frac{\left|  h\right|  ^{2}}{\left(  R+\rho\right)  ^{2}%
}\right) \label{bochner}\\
&  +4C\rho\frac{\left|  h\right|  ^{2}}{\left(  R+\rho\right)  ^{2}}.\nonumber
\end{align}
\emph{Consequently, by direct application of the maximum principle, the norm
of the solution to the linearized Ricci flow equation is comparable to the
scalar curvature plus a constant: }%
\[
\left|  h\right|  \leq C\left(  R+\rho\right)
\]
\emph{on }$M\times\lbrack0,T),$\emph{ where }$C$\emph{ depends only on
}$g\left(  0\right)  ,$\emph{ }$\rho$\emph{ and }$T.\medskip$
\emph{Furthermore, when }$\rho=0,$ $C$\emph{ is independent of }$T.\medskip$

Taking $\rho=0$ and $h_{ij}=R_{ij}\,,$ we obtain:

\begin{corollary}
(Hamilton, \cite{H1}) If $\left(  M^{3},g\left(  t\right)  \right)  ,$
$t\in\lbrack0,T),$ $T<\infty$ is a solution to the Ricci flow on a closed
$3$-manifold with positive scalar curvature, then there exists a constant
$C<\infty$ such that%
\[
\frac{\left|  Rc\right|  }{R}\leq C
\]
on $M\times\lbrack0,T).$
\end{corollary}

In \S 2 we recall how the system (\ref{ricciflow})-(\ref{LichLapEq}) is
obtained by linearizing the Ricci flow using a version of DeTurck's trick with
a time-dependent background metric. In \S 3 we give the proof of equation
(\ref{bochner}), from which the main theorem follows. This depends on the
nonnegativity of a certain degree 4 polynomial in 6 variables (Lemma
\ref{quartic nonneg}), which is proved in \S 4.

\section{The linearized Ricci flow system\label{defmot}}

This section is mainly to motivate our study of the linearized Ricci flow
system. The reader well familiar with DeTurck's trick \cite{D} may skip this
section. 

As we stated in the introduction, a solution to the\emph{ linearized Ricci
flow system} consists of a complete solution $\left(  M^{n},g_{o}\left(
t\right)  \right)  ,$ $t\in\lbrack0,T),$ to the Ricci flow%
\begin{equation}
\frac{\partial}{\partial t}g_{ij}=-2R_{ij}\label{RF}%
\end{equation}
coupled with a solution $h\left(  t\right)  ,$ $t\in\lbrack0,T),$ to the
\emph{Lichnerowicz laplacian heat equation}
\begin{equation}
\frac{\partial}{\partial t}h_{ij}=\left(  \Delta_{L}h\right)  _{ij}%
\label{LLHE}%
\end{equation}
where%
\[
\left(  \Delta_{L}h\right)  _{ij}\doteqdot\Delta h_{ij}+2R_{kijl}h_{kl}%
-R_{ik}h_{kj}-R_{jk}h_{ki}.
\]
The Lichnerowicz laplacian $\Delta_{L}$ is defined using the evolving metric
$g_{o}\left(  t\right)  .$ Our main interest is when the solution is compact.
However, in view of compactness arguments in the category of pointed
solutions, it may be of interest to study the linearized Ricci flow system for
complete, noncompact solutions.

This system arises as follows. Given a solution $\left(  M^{n},g_{o}\left(
t\right)  \right)  ,$ $t\in\lbrack0,T_{o}),$ to the Ricci flow, consider the
\emph{modified Ricci flow with time-dependent background metrics}
$g_{o}\left(  t\right)  \,$:%
\begin{equation}
\frac{\partial}{\partial t}g_{ij}=-2R_{ij}+\nabla_{i}W_{j}+\nabla_{j}W_{i}
\label{MRF-1E}%
\end{equation}
where the 1-forms $W\left(  t\right)  $ are defined by%
\begin{equation}
W\left(  t\right)  _{\ell}\doteqdot g\left(  t\right)  _{\ell k}g\left(
t\right)  ^{ij}\left(  \Gamma\left[  g\left(  t\right)  \right]  _{ij}%
^{k}-\Gamma\left[  g_{o}\left(  t\right)  \right]  _{ij}^{k}\right)
\label{W1}%
\end{equation}
and the covariant derivatives are with respect to the metrics $g\left(
t\right)  .$ This is DeTurck's trick \cite{D} with a \emph{time-dependent}
background metric. Note that the solution $g_{o}\left(  t\right)  $ to the
Ricci flow is itself also a solution the modified Ricci flow with background
metrics $g_{o}\left(  t\right)  .$

There exists a unique solution to the initial value problem for the modified
Ricci flow for short time. The \emph{modified Ricci tensor}, which we define
to be the rhs of (\ref{MRF-1E}), depends only on $g\left(  t\right)  $ and
$g_{o}\left(  t\right)  .$ We first compute the linearization of the modified
Ricci flow about the solution $g_{o}\left(  t\right)  .$ In this case we get
the Lichnerowicz laplacian heat equation. Hence the modified Ricci flow is a
parabolic equation, which in turn, implies uniqueness and short time
existence. DeTurck gave this argument as a new proof of the short time
existence and uniqueness of solutions to the Ricci flow originally proved in
\cite{H1}.

Let $\left\{  g_{s,o}\right\}  _{s\in\left(  -\varepsilon,\varepsilon\right)
}$ be a smooth, one-parameter family of initial metrics with $g_{0,o}%
=g_{o}\left(  0\right)  .$ Consider the one-parameter family $\left\{
g_{s}\left(  t\right)  ,t\in\lbrack0,T_{s})\right\}  _{s\in\left(
-\varepsilon,\varepsilon\right)  }$ of solutions to the modified Ricci flow:
\begin{align}
\frac{\partial}{\partial t}\left(  g_{s}\right)  _{ij}  &  =-2\left(
Rc_{s}\right)  _{ij}+\nabla_{i}W_{j}+\nabla_{j}W_{i}\label{MRF-E}\\
g_{s}\left(  0\right)   &  =g_{s,o}, \label{MRF-IV}%
\end{align}
where
\begin{equation}
W_{s}\left(  t\right)  _{\ell}\doteqdot g_{s}\left(  t\right)  _{\ell k}%
g_{s}\left(  t\right)  ^{ij}\left(  \Gamma\left[  g_{s}\left(  t\right)
\right]  _{ij}^{k}-\Gamma\left[  g_{o}\left(  t\right)  \right]  _{ij}%
^{k}\right)  , \label{W}%
\end{equation}
and the Ricci tensor and covariant derivative are with respect to the metrics
$g_{s}\left(  t\right)  .$ Recall that
\[
g_{0}\left(  t\right)  \equiv g_{o}\left(  t\right)  \quad\text{and\quad}%
W_{0}\left(  t\right)  \equiv0.
\]

Define%
\begin{equation}
v_{ij}\left(  t\right)  \doteqdot\left.  \frac{\partial}{\partial s}\right|
_{s=0}g_{s}\left(  t\right)  _{ij}. \label{vij}%
\end{equation}
We shall call $v$ the \emph{variation of the metric tensor}.

Let $V\left(  t\right)  =g_{0}\left(  t\right)  ^{ij}v\left(  t\right)
_{ij}.$ A standard computation yields (see for example \cite{H1})

\begin{lemma}
The \emph{variation of the Ricci tensor} is%
\begin{align}
&  \left.  \frac{\partial}{\partial s}\right|  _{s=0}\left(  -2Rc\left[
g_{s}\left(  t\right)  \right]  _{ij}\right) \nonumber\\
&  =\Delta v_{ij}+2R_{kijl}v_{kl}-R_{ik}v_{kj}-R_{jk}v_{ki}+\nabla_{i}%
\nabla_{j}V-\nabla_{i}\nabla^{k}v_{kj}-\nabla_{j}\nabla^{k}v_{ki}\nonumber\\
&  =\left(  \Delta_{L}v\right)  _{ij}+\nabla_{i}\nabla_{j}V-\nabla_{i}%
\nabla^{k}v_{kj}-\nabla_{j}\nabla^{k}v_{ki} \label{VRT}%
\end{align}

\end{lemma}

Recall the algebraic \emph{Einstein operator}
\[
G\left(  v\right)  _{ij}\doteqdot v_{ij}-\frac{1}{2}Vg_{ij}%
\]
(which takes the Ricci tensor to the Einstein tensor: $G\left(  Rc\right)
_{ij}=R_{ij}-\frac{1}{2}Rg_{ij}$) and the \emph{divergence}%
\[
\delta:C^{\infty}\left(  S^{2}T^{\ast}M\right)  \rightarrow C^{\infty}\left(
T^{\ast}M\right)
\]
where%
\[
\delta\left(  T\right)  _{i}\doteqdot-g^{jk}\nabla_{k}T_{ij}.
\]
The $L^{2}$ \emph{adjoint} of $\delta$%
\[
\delta^{\ast}:C^{\infty}\left(  T^{\ast}M\right)  \rightarrow C^{\infty
}\left(  S^{2}T^{\ast}M\right)
\]
is the same as the \emph{Lie derivative operator} acting on the metric:%
\[
\delta^{\ast}\left(  v\right)  _{ij}=\frac{1}{2}\left(  \nabla_{i}v_{j}%
+\nabla_{j}v_{i}\right)  =\frac{1}{2}\left(  L_{v^{\ast}}g\right)  _{ij}%
\]
where $\left(  v^{\ast}\right)  ^{i}=g^{ij}v_{j}.$ The last three terms on the
rhs of equation (\ref{VRT}) may be rewritten as%
\[
\nabla_{i}\nabla_{j}V-\nabla_{i}\nabla^{k}v_{kj}-\nabla_{j}\nabla^{k}%
v_{ki}=2\left[  \delta^{\ast}\left(  \delta\left[  G\left(  v\right)  \right]
\right)  \right]  _{ij}.
\]
Hence we have the following well-known identity:

\begin{lemma}
The variation of the Ricci tensor has the form%
\[
\left.  \frac{\partial}{\partial s}\right|  _{s=0}\left(  -2Rc\left[
g_{s}\left(  t\right)  \right]  _{ij}\right)  =\left(  \Delta_{L}v\right)
_{ij}+2\left[  \delta^{\ast}\left(  \delta\left[  G\left(  v\right)  \right]
\right)  \right]  _{ij}.
\]

\end{lemma}

We may rewrite the 1-form $W_{s}\left(  t\right)  $ as
\begin{align*}
W_{s}\left(  t\right)  _{\ell}  &  =-\frac{1}{2}g_{s}\left(  t\right)  _{\ell
k}g_{0}\left(  t\right)  ^{kp}g_{s}\left(  t\right)  ^{ij}\left(  \nabla
_{i}\left[  g_{0}\right]  _{jp}+\nabla_{j}[g_{0}]_{ip}-\nabla_{p}\left[
g_{0}\right]  _{ij}\right) \\
&  =W_{s}\left(  t\right)  _{\ell}=\frac{1}{2}g_{s}\left(  t\right)  _{\ell
k}g_{0}\left(  t\right)  ^{kp}\left(  \delta\left[  G\left(  g_{0}\right)
\right]  \right)  _{p}\,,
\end{align*}
where the covariant derivatives are with respect to the metrics $g_{s}\left(
t\right)  .$ Define%
\[
\left(  g_{0}\right)  ^{-1}:C^{\infty}\left(  T^{\ast}M\right)  \rightarrow
C^{\infty}\left(  T^{\ast}M\right)
\]
by%
\[
\left(  g_{0}\right)  ^{-1}T\doteqdot g_{s}\left(  t\right)  _{\ell k}%
g_{0}\left(  t\right)  ^{kp}T_{p}.
\]
Then%
\[
W_{s}\left(  t\right)  _{\ell}=\frac{1}{2}\left[  \left(  g_{0}\right)
^{-1}\left(  \delta\left[  G\left(  g_{0}\right)  \right]  \right)  \right]
_{\ell}.
\]
Hence the last two terms of the modified Ricci tensor can be expressed as%
\[
\nabla_{i}W_{j}+\nabla_{j}W_{i}=2\left(  \delta^{\ast}W\right)  _{ij}=\left(
\delta^{\ast}\left[  \left(  g_{0}\right)  ^{-1}\left(  \delta\left[  G\left(
g_{0}\right)  \right]  \right)  \right]  \right)  _{ij},
\]
so that%
\[
-2R_{ij}+\nabla_{i}W_{j}+\nabla_{j}W_{i}=-2R_{ij}+\left(  \delta^{\ast}\left[
\left(  g_{0}\right)  ^{-1}\left(  \delta\left[  G\left(  g_{0}\right)
\right]  \right)  \right]  \right)  _{ij}.
\]

Thus the only change to DeTurck's modification of the Ricci flow that we have
made is that we allow the background metric $g_{0}$ to depend on time. In
particular, we take $g_{0}$ to be the solution of the Ricci flow that we are
linearizing about.

The motivation for studying the linearized Ricci flow system is the following.

\begin{proposition}
The variation $v\left(  t\right)  $ of the metric tensor $g\left(  t\right)  $
corresponding to a one-parameter family $\left\{  g_{s}\left(  t\right)
\right\}  _{s\in\left(  -\varepsilon,\varepsilon\right)  }$ of solutions to
the modified Ricci flow (\ref{MRF-E})-(\ref{MRF-IV}) is a solution to the
Lichnerowicz laplacian heat equation:%
\[
\frac{\partial}{\partial t}v=\Delta_{L}v.
\]
That is, the pair $\left(  g,v\right)  $ is a solution to the linearized Ricci
flow system (\ref{RF})-(\ref{LLHE}).
\end{proposition}

\textbf{Proof}. We compute%
\begin{align*}
\frac{\partial}{\partial t}v_{ij}\left(  t\right)   &  =\frac{\partial
}{\partial t}\left.  \frac{\partial}{\partial s}\right|  _{s=0}g_{s}\left(
t\right)  _{ij}=\left.  \frac{\partial}{\partial s}\right|  _{s=0}%
\frac{\partial}{\partial t}g_{s}\left(  t\right)  _{ij}\\
&  =\left.  \frac{\partial}{\partial s}\right|  _{s=0}\left(  -2R_{ij}%
+\nabla_{i}W_{j}+\nabla_{j}W_{i}\right) \\
&  =\left.  \frac{\partial}{\partial s}\right|  _{s=0}\left(  -2Rc\left[
g_{s}\left(  t\right)  \right]  _{ij}\right)  +\nabla_{i}\left(  \left.
\frac{\partial}{\partial s}\right|  _{s=0}W_{s}\left(  t\right)  _{j}\right)
+\nabla_{j}\left(  \left.  \frac{\partial}{\partial s}\right|  _{s=0}%
W_{s}\left(  t\right)  _{i}\right) \\
&  =\left(  \Delta_{L}v\right)  _{ij}+\nabla_{i}\nabla_{j}V-\nabla_{i}%
\nabla^{k}v_{kj}-\nabla_{j}\nabla^{k}v_{ki}\\
&  +\nabla_{i}\left(  \frac{1}{2}g_{0}\left(  t\right)  ^{k\ell}\left(
\nabla_{k}v_{\ell j}+\nabla_{\ell}v_{kj}-\nabla_{j}v_{k\ell}\right)  \right)
\\
&  +\nabla_{j}\left(  \frac{1}{2}g_{0}\left(  t\right)  ^{k\ell}\left(
\nabla_{k}v_{\ell i}+\nabla_{\ell}v_{kj}-\nabla_{i}v_{k\ell}\right)  \right)
\\
&  =\left(  \Delta_{L}v\right)  _{ij},
\end{align*}
where we used $W_{0}\left(  t\right)  \equiv0$ to obtain the fourth equality,
and the fact that $\nabla_{i}\nabla_{j}V=\nabla_{j}\nabla_{i}V$ for the last equality.

\section{Proof of the pinching estimate}

In this and the following section, we give the derivation of equation
(\ref{bochner}), which implies the main theorem. Using the Ricci flow
equation, (\ref{LichLapEq}), and the standard equation $\frac{\partial
}{\partial t}R=\Delta R+2\left|  Rc\right|  ^{2},$ we compute%
\begin{align*}
\frac{\partial}{\partial t}\left(  \frac{\left|  h\right|  ^{2}}{\left(
R+\rho\right)  ^{2}}\right)   &  =-\frac{2}{\left(  R+\rho\right)  ^{2}}%
\frac{\partial}{\partial t}g_{ij}\cdot h_{ij}^{2}+\frac{2}{\left(
R+\rho\right)  ^{2}}\left(  \frac{\partial}{\partial t}h\right)  \cdot
h-2\frac{\left|  h\right|  ^{2}}{\left(  R+\rho\right)  ^{3}}\frac{\partial
}{\partial t}R\\
&  =\frac{4Rc\cdot h^{2}}{\left(  R+\rho\right)  ^{2}}+\frac{2}{\left(
R+\rho\right)  ^{2}}\left(  \Delta h_{ij}+2R_{kijl}h_{kl}-R_{ik}h_{kj}%
-R_{jk}h_{ki}\right)  h_{ij}\\
&  -2\frac{\left|  h\right|  ^{2}}{\left(  R+\rho\right)  ^{3}}\left(  \Delta
R+2\left|  Rc\right|  ^{2}\right) \\
&  =\Delta\left(  \frac{\left|  h\right|  ^{2}}{\left(  R+\rho\right)  ^{2}%
}\right)  -2\frac{\left|  \nabla h\right|  ^{2}}{\left(  R+\rho\right)  ^{2}%
}-6\frac{\left|  h\right|  ^{2}}{\left(  R+\rho\right)  ^{4}}\left|  \nabla
R\right|  ^{2}\\
&  +\frac{8}{\left(  R+\rho\right)  ^{3}}h\nabla_{i}R\cdot\nabla_{i}h+\frac
{4}{\left(  R+\rho\right)  ^{2}}R_{ijkl}h_{il}h_{jk}-\frac{4}{\left(
R+\rho\right)  ^{3}}\left|  h\right|  ^{2}\left|  Rc\right|  ^{2}.
\end{align*}
Since%
\[
\frac{1}{R+\rho}\nabla R\cdot\nabla\left(  \frac{\left|  h\right|  ^{2}%
}{\left(  R+\rho\right)  ^{2}}\right)  =-2\frac{\left|  h\right|  ^{2}%
}{\left(  R+\rho\right)  ^{4}}\left|  \nabla R\right|  ^{2}+\frac{2}{\left(
R+\rho\right)  ^{3}}h\nabla_{i}R\cdot\nabla_{i}h,
\]
we may rewrite the above evolution equation as%
\begin{align}
\frac{\partial}{\partial t}\left(  \frac{\left|  h\right|  ^{2}}{\left(
R+\rho\right)  ^{2}}\right)   &  =\Delta\left(  \frac{\left|  h\right|  ^{2}%
}{\left(  R+\rho\right)  ^{2}}\right)  +\frac{2}{R+\rho}\nabla R\cdot
\nabla\left(  \frac{\left|  h\right|  ^{2}}{\left(  R+\rho\right)  ^{2}%
}\right) \label{norm h over R}\\
&  -2\frac{\left|  \left(  R+\rho\right)  \nabla_{i}h_{jk}-\nabla_{i}%
Rh_{jk}\right|  ^{2}}{\left(  R+\rho\right)  ^{4}}+\frac{4}{\left(
R+\rho\right)  ^{2}}R_{ijkl}h_{il}h_{jk}\nonumber\\
&  -\frac{4}{\left(  R+\rho\right)  ^{3}}\left|  h\right|  ^{2}\left|
Rc\right|  ^{2}.\nonumber
\end{align}

When $n=3,$ we have the identity%
\[
R_{ijkl}=R_{il}g_{jk}+R_{jk}g_{il}-R_{ik}g_{jl}-R_{jl}g_{ik}-\frac{1}%
{2}R\left(  g_{il\,}g_{jk\,-}g_{ik\,}g_{jl}\right)  .
\]
Hence%
\begin{align*}
\frac{\partial}{\partial t}\left(  \frac{\left|  h\right|  ^{2}}{\left(
R+\rho\right)  ^{2}}\right)   &  =\Delta\left(  \frac{\left|  h\right|  ^{2}%
}{\left(  R+\rho\right)  ^{2}}\right)  +\frac{2}{R+\rho}\nabla R\cdot
\nabla\left(  \frac{\left|  h\right|  ^{2}}{\left(  R+\rho\right)  ^{2}%
}\right) \\
&  -\frac{2}{\left(  R+\rho\right)  ^{4}}\left|  \left(  R+\rho\right)
\nabla_{i}h_{jk}-\nabla_{i}Rh_{jk}\right|  ^{2}+4P
\end{align*}
where
\begin{align*}
P  &  =\frac{1}{\left(  R+\rho\right)  ^{2}}\left[  R_{il}g_{jk}+R_{jk}%
g_{il}-R_{ik}g_{jl}-R_{jl}g_{ik}-\frac{1}{2}R\left(  g_{il\,}g_{jk\,-}%
g_{ik\,}g_{jl}\right)  \right]  h_{il}h_{jk}\\
&  -\frac{1}{\left(  R+\rho\right)  ^{3}}\left|  h\right|  ^{2}\left|
Rc\right|  ^{2}\\
&  =\frac{1}{\left(  R+\rho\right)  ^{3}}\left[  2\left(  R+\rho\right)
Rc\cdot hH-2\left(  R+\rho\right)  Rc\cdot h^{2}+\frac{R}{2}\left(
R+\rho\right)  \left(  \left|  h\right|  ^{2}-H^{2}\right)  -\left|  h\right|
^{2}\left|  Rc\right|  ^{2}\right] \\
&  =\frac{1}{\left(  R+\rho\right)  ^{3}}\left[  2RRc\cdot hH-2RRc\cdot
h^{2}+\frac{1}{2}R^{2}\left(  \left|  h\right|  ^{2}-H^{2}\right)  -\left|
h\right|  ^{2}\left|  Rc\right|  ^{2}\right] \\
&  +\frac{\rho}{\left(  R+\rho\right)  ^{3}}\left[  2Rc\cdot hH-2Rc\cdot
h^{2}+\frac{1}{2}R\left(  \left|  h\right|  ^{2}-H^{2}\right)  \right]  .
\end{align*}

The main theorem is now a consequence of the following inequality, which we
shall prove in the next section.

\begin{lemma}
\label{quartic nonneg}We have for any metric $g$ and symmetric $2$-tensor $h,$
the inequality%
\[
\left|  h\right|  ^{2}\left|  Rc\right|  ^{2}-2RHRc\cdot h+2RRc\cdot
h^{2}+\frac{1}{2}R^{2}\left(  H^{2}-\left|  h\right|  ^{2}\right)  \geq0.
\]

\end{lemma}

\noindent This is because, then%
\begin{align*}
\frac{\partial}{\partial t}\left(  \frac{\left|  h\right|  ^{2}}{\left(
R+\rho\right)  ^{2}}\right)   &  \leq\Delta\left(  \frac{\left|  h\right|
^{2}}{\left(  R+\rho\right)  ^{2}}\right)  +\frac{2}{R+\rho}\nabla
R\cdot\nabla\left(  \frac{\left|  h\right|  ^{2}}{\left(  R+\rho\right)  ^{2}%
}\right) \\
&  +4\frac{\rho}{\left(  R+\rho\right)  ^{3}}\left[  2Rc\cdot hH-2Rc\cdot
h^{2}+\frac{1}{2}R\left(  \left|  h\right|  ^{2}-H^{2}\right)  \right]  .
\end{align*}
On the other hand, we have the estimate $\left|  Rc\right|  \leq C\left(
R+\rho\right)  $ (see Theorem 24.4 of \cite{H2}), which implies%
\[
\frac{\rho}{\left(  R+\rho\right)  ^{3}}\left[  2Rc\cdot hH-2Rc\cdot
h^{2}+\frac{1}{2}R\left(  \left|  h\right|  ^{2}-H^{2}\right)  \right]  \leq
C\rho\frac{\left|  h\right|  ^{2}}{\left(  R+\rho\right)  ^{2}}.
\]
Hence%
\[
\frac{\partial}{\partial t}\left(  \frac{\left|  h\right|  ^{2}}{\left(
R+\rho\right)  ^{2}}\right)  \leq\Delta\left(  \frac{\left|  h\right|  ^{2}%
}{\left(  R+\rho\right)  ^{2}}\right)  +\frac{2}{R+\rho}\nabla R\cdot
\nabla\left(  \frac{\left|  h\right|  ^{2}}{\left(  R+\rho\right)  ^{2}%
}\right)  +4C\rho\frac{\left|  h\right|  ^{2}}{\left(  R+\rho\right)  ^{2}}.
\]
If $t\in\lbrack0,T),$ then applying the maximum principle implies%
\[
\frac{\left|  h\right|  ^{2}}{\left(  R+\rho\right)  ^{2}}\left(  t\right)
\leq C_{0}\exp\left(  4C\rho T\right)  ,
\]
where $C_{0}=\max_{t=0}\left|  h\right|  ^{2}/\left(  R+\rho\right)  ^{2}.$
\textbf{q.e.d.}

\section{Nonnegativity of a degree $4$ homogeneous polynomial in $6$
variables}

\noindent\textbf{Proof of Lemma \ref{quartic nonneg}}. Since $h$ is symmetric,
we may assume $h$ is diagonal. Let $h_{1},h_{2},h_{3}$ denote the eigenvalues
of $h$ and let $r_{1}=R_{11},r_{2}=R_{22},r_{3}=R_{33}$ denote the diagonal
entries of $R_{ij}.$ Then%
\begin{align*}
-R^{3}P  &  =\left|  h\right|  ^{2}\left|  Rc\right|  ^{2}-2RHRc\cdot
h+2RRc\cdot h^{2}+\frac{1}{2}R^{2}\left(  H^{2}-\left|  h\right|  ^{2}\right)
\\
&  \geq Q\doteqdot\left(  h_{1}^{2}+h_{2}^{2}+h_{3}^{2}\right)  \left(
r_{1}^{2}+r_{2}^{2}+r_{3}^{2}\right) \\
&  -2\left(  r_{1}+r_{2}+r_{3}\right)  \left(  h_{1}+h_{2}+h_{3}\right)
\left(  r_{1}h_{1}+r_{2}h_{2}+r_{3}h_{3}\right) \\
&  +2\left(  r_{1}+r_{2}+r_{3}\right)  \left(  r_{1}h_{1}^{2}+r_{2}h_{2}%
^{2}+r_{3}h_{3}^{2}\right)  +\left(  r_{1}+r_{2}+r_{3}\right)  ^{2}\left(
h_{1}h_{2}+h_{1}h_{3}+h_{2}h_{3}\right)  ,
\end{align*}
where we used the inequality (throwing away the off-diagonal entries of
$R_{ij}$)%
\[
\left|  Rc\right|  ^{2}\geq r_{1}^{2}+r_{2}^{2}+r_{3}^{2}.
\]
We expand and simplify this as%
\begin{align*}
&  Q=r_{1}^{2}h_{1}^{2}+r_{2}^{2}h_{2}^{2}+r_{3}^{2}h_{3}^{2}+r_{1}^{2}%
h_{2}^{2}+r_{1}^{2}h_{3}^{2}+r_{2}^{2}h_{1}^{2}+r_{2}^{2}h_{3}^{2}+r_{3}%
^{2}h_{1}^{2}+r_{3}^{2}h_{2}^{2}\\
&  +r_{3}^{2}h_{1}h_{2}+r_{1}^{2}h_{2}h_{3}+r_{2}^{2}h_{1}h_{3}\\
&  -r_{1}^{2}h_{1}h_{2}-r_{1}^{2}h_{1}h_{3}-r_{2}^{2}h_{1}h_{2}-r_{2}^{2}%
h_{2}h_{3}-r_{3}^{2}h_{1}h_{3}-r_{3}^{2}h_{2}h_{3}\\
&  -2r_{1}r_{2}h_{1}h_{2}-2r_{1}r_{3}h_{1}h_{3}-2r_{2}r_{3}h_{2}h_{3}.
\end{align*}
Writing $Q$ as a bilinear form in $h=\left(  h_{1},h_{2},h_{3}\right)  $ with
coefficients in $r=\left(  r_{1},r_{2},r_{3}\right)  ,$ we have%
\begin{align*}
&  Q=\left(  r_{1}^{2}+r_{2}^{2}+r_{3}^{2}\right)  h_{1}^{2}+\left(  r_{1}%
^{2}+r_{2}^{2}+r_{3}^{2}\right)  h_{2}^{2}+\left(  r_{1}^{2}+r_{2}^{2}%
+r_{3}^{2}\right)  h_{3}^{2}\\
&  +\left(  -r_{1}^{2}-r_{2}^{2}+r_{3}^{2}-2r_{1}r_{2}\right)  h_{1}%
h_{2}+\left(  r_{1}^{2}-r_{2}^{2}-r_{3}^{2}-2r_{2}r_{3}\right)  h_{2}h_{3}\\
&  +\left(  -r_{1}^{2}+r_{2}^{2}-r_{3}^{2}-2r_{1}r_{3}\right)  h_{1}h_{3}.
\end{align*}
The lemma follows from the claim that the polynomial $Q$ takes nonnegative
values for all real values of $r_{1},r_{2},r_{3},h_{1},h_{2},h_{3}$. To prove
the claim, it is enough to show that the Hessian matrix
\[%
\begin{array}
[c]{rcl}%
H & := & \displaystyle\left[  \frac{\partial^{2}Q}{\partial h_{i}\partial
h_{j}}\right] \\
&  & \\
& = & \left[
\begin{array}
[c]{ccc}%
2r_{1}^{2}+2r_{2}^{2}+2r_{3}^{2} & r_{3}^{2}-r_{1}^{2}-r_{2}^{2}-2\,r_{{1}%
}r_{{2}} & r_{2}^{2}-r_{1}^{2}-r_{3}^{2}-2\,r_{{1}}r_{{3}}\\
r_{3}^{2}-r_{1}^{2}-r_{2}^{2}-2\,r_{{1}}r_{{2}} & 2r_{1}^{2}+2r_{2}^{2}%
+2r_{3}^{2} & r_{1}^{2}-r_{2}^{2}-r_{3}^{2}-2\,r_{{2}}r_{{3}}\\
r_{2}^{2}-r_{1}^{2}-r_{3}^{2}-2\,r_{{1}}r_{{3}} & r_{1}^{2}-r_{2}^{2}%
-r_{3}^{2}-2\,r_{{2}}r_{{3}} & 2r_{1}^{2}+2r_{2}^{2}+2r_{3}^{2}%
\end{array}
\right]
\end{array}
\]
is positive semidefinite for all real values of $r_{1},r_{2},r_{3}$. In turn,
it is enough to show that the determinants
\[
\Delta_{1}:=H_{11},\;\;\;\Delta_{2}:=\left|
\begin{array}
[c]{cc}%
H_{11} & H_{12}\\
H_{21} & H_{22}%
\end{array}
\right|  ,\;\;\;\Delta_{3}:=\left|
\begin{array}
[c]{ccc}%
H_{11} & H_{12} & H_{13}\\
H_{21} & H_{22} & H_{23}\\
H_{31} & H_{32} & H_{33}%
\end{array}
\right|
\]
take nonnegative values for all real values of $r_{1},r_{2},r_{3}$. The
nonnegativity of
\[
\Delta_{1}=2r_{1}^{2}+2r_{2}^{2}+2r_{3}^{2}%
\]
is clear. The nonnegativity of
\[
\Delta_{2}=(2r_{1}+2r_{2}+2r_{3})^{2}-(r_{1}^{2}-r_{2}^{2}-2r_{2}r_{3}%
-r_{3}^{2})^{2}%
\]
follows from the inequality
\[
2r_{1}^{2}+2r_{2}^{2}+2r_{3}^{2}\geq r_{1}^{2}+r_{2}^{2}+2|r_{1}r_{2}%
|+r_{3}^{2}\geq|r_{1}^{2}-r_{2}^{2}-2r_{2}r_{3}-r_{3}^{2}|.
\]
One has
\[
3\Delta_{3}=2\,{X}^{6}-6\,{X}^{4}Y-24\,{X}^{2}{Y}^{2}+24\,{Y}^{3}+16\,{X}%
^{3}Z,
\]
where
\[
X:=r_{1}+r_{2}+r_{3},\;\;\;Y:=r_{1}^{2}+r_{2}^{2}+r_{3}^{2},\;\;\;Z:=r_{1}%
^{3}+r_{2}^{3}+r_{3}^{3}.
\]
The reader will have no difficulty verifying this identity with the help of a
computer algebra package. Now on any circle in $(r_{1},r_{2},r_{3})$-space
defined by holding $X$ and $Y$ constant, the extremal values of $Z$ occur at
points with two coordinates equal, as one verifies by a straightforward
Lagrange multiplier argument. Since for fixed $X$ and $Y$, $\Delta_{3}$
depends linearly on $Z$, it follows that on any circle defined by holding $X$
and $Y$ constant, the extreme values of $\Delta_{3}$ are taken at points with
two coordinates equal. After making the evident reductions, the inequality
\[
\Delta_{3}(x,x,1)=4(8x^{2}+1)(x-1)^{2}\geq0
\]
suffices to prove the nonnegativity of $\Delta_{3}$, and in turn the claim.
With hindsight, the reader can see that it was overkill to actually write out
the expression of $\Delta_{3}$ in terms of $X$, $Y$ and $Z$; all that was used
in the proof of nonnegativity of $\Delta_{3}$ was the fact that for fixed $X$
and $Y$, $\Delta_{3}$ depends linearly on $Z$. \textbf{q.e.d.\medskip}

\textbf{Acknowledgment. }We would like to thank Jiaping Wang for bringing to
our attention Matt Gursky's work which partly inspired the result in this paper.

\end{document}